\documentclass[12pt]{amsart}
\usepackage{mathrsfs}
\usepackage{amsfonts}
\usepackage{amssymb}
\usepackage{amsfonts, amscd, amsmath, mathrsfs, amssymb, amsthm, amsxtra, bbding, epsfig, graphicx, latexsym, url, mathbbol, bbold}
\usepackage[papersize={8in,11in},textwidth=15.3cm,textheight=21cm,centering]{geometry}
\usepackage{enumerate}

\usepackage{xcolor}
\definecolor{cite}{rgb}{0.00,0.60,1.00}
\definecolor{url}{rgb}{1.00,0.10,0.80}
\definecolor{link}{rgb}{0.00,0.00,1.00}
\usepackage[colorlinks,linkcolor=link,urlcolor=url,citecolor=cite,pagebackref,breaklinks]{hyperref}

\hypersetup{
pdfstartpage=1,
pdfstartview=FitH}

\DeclareFontFamily{U}{mathx}{\hyphenchar\font45}
\DeclareFontShape{U}{mathx}{m}{n}{
      <5> <6> <7> <8> <9> <10>
      <10.95> <12> <14.4> <17.28> <20.74> <24.88>
      mathx10
      }{}
\DeclareSymbolFont{mathx}{U}{mathx}{m}{n}
\DeclareMathAccent{\widecheck}{\mathalpha}{mathx}{"71}

 \usepackage{caption} 
\numberwithin{equation}{section}

\allowdisplaybreaks

\newtheorem{lemma}{Lemma}[section]

\makeatletter
\newcounter{roem}
\renewcommand{\theroem}{\Roman{roem}}

\newcommand{\c@org@eq}{}
\let\c@org@eq\c@equation
\newcommand{\org@theeq}{}
\let\org@theeq\theequation

\newcommand{\setroem}{
\let\c@equation\c@roem
 \let\theequation\theroem}

\newcommand{\setarab}{
\let\c@equation\c@org@eq
\let\theequation\org@theeq}
\makeatother

\newtheorem*{claim*}{Claim}

\theoremstyle{remark}

\newcommand{\ue}{\mathrm{e}}

\DeclareMathOperator{\Mod}{mod}

\renewcommand{\bmod}[1]{\,(\Mod{ #1})}

\newcommand{\bZ}{\mathbf{Z}}

\usepackage{graphicx}
\usepackage{tikz}

\begin{document}

\title{An elementary proof of the Selberg identity for Kloosterman sums}

\author{Ping Xi}

\address{School of Mathematics and Statistics, Xi'an Jiaotong University, Xi'an 710049, P. R. China}

\email{ping.xi@xjtu.edu.cn}

\begin{abstract} 
We give an elementary proof of the Selberg identity for Kloosterman sums, which only requires the orthogonality of additive characters.
\end{abstract}

\maketitle

\setcounter{tocdepth}{1}

\section{Introduction}\label{sec:Introduction}

Let us recall the classical Selberg identity
\begin{align}\label{eq:KSidentity}
S(m,n;c)=\sum_{d|(m,n,c)}d S(mnd^{-2},1;cd^{-1}).
\end{align}
Here $S(m,n;c)$ is the Kloosterman sum defined by
\begin{align*}
S(m,n;c)=\sideset{}{^*}\sum_{a\bmod c}\ue\Big(\frac{ma+n\overline{a}}{c}\Big)
\end{align*}
with $c\geqslant1,$ $m,n\in\bZ$ and $a\overline{a}\equiv1
\bmod c.$
The identity can be regarded as a reflection of the multiplicativitiy of Hecke operators, and it was announced by Selberg \cite{Se38} early in 1930's. The published proof was not found until Kuznestsov \cite{Ku80} developed his trace formula in late 1970's, from which \eqref{eq:KSidentity} was concluded. Namely, the first proof of \eqref{eq:KSidentity} utilizes automorphic forms, and it is natural to ask if there exists a non-automorphic proof.
In particular, elementary proofs are highly desired. Of course, we are not valid here to give an appropriate definition of ``elementary"; this should depend on the real nature of Kloosterman sums, which seems to be still mysterious at present.

Suppose we now have some intuition on ``elementary", and there exist at least three elementary proofs of \eqref{eq:KSidentity}:

\begin{itemize}
\item  Matthes \cite{Ma90} gave a proof using the Chinese remainder theorem and M\"obius inversion.
\item Andersson \cite{An06} gave another proof, in which the theory of Hecke operators was applied. 
\item Harcos and K\'arolyi \cite{HK17} proved \eqref{eq:KSidentity} by direct computations for Kloosterman sums to prime power moduli, before which they also apply the 
Chinese remainder theorem.

\end{itemize}

In this note we give another elementary and direct proof, which only requires the orthogonality of additive characters. The argument here is not new, in fact it has been included as an appendix to \cite{Xi18}, in which the primary aim is to derive a similar identity for hyper-Kloosterman sums following the similar spirit. We now decide to make it as a separate material since it might be of independent interests.

\section{The proof}
The following lemma contributes as the only tool in the subsequent arguments.

\begin{lemma}\label{lm:orthogonality}
Let $q$ be a positive integer. For each $n\in\bZ,$ we have
\begin{align*}
\frac{1}{q}\sum_{a\bmod q}\ue\Big(\frac{na}{q}\Big)
=\begin{cases}
1,\ \ & q\mid n,\\
0,\ \ & q\nmid n.
\end{cases}
\end{align*}

\end{lemma}

We now start the proof. From Lemma \ref{lm:orthogonality} it follows that
\begin{align*}
S(m,n;c)
&=\frac{1}{c}\sum_{a\bmod c}\mathop{\sum\sum}_{x,y\bmod c}\ue\Big(\frac{x(m+ay)+ny-a}{c}\Big)\\
&=\frac{1}{c}\sum_{d\mid c}\sum_{\substack{a\bmod{c/d}\\(a,c/d)=1}}\mathop{\sum\sum}_{x,y\bmod c}\ue\Big(\frac{x(m+ady)+ny-ad}{c}\Big).\end{align*}
We evaluate the $x$-sum by orthogonality, getting
\begin{align*}
S(m,n;c)&=\sum_{d\mid c}\sum_{\substack{a\bmod{c/d}\\(a,c/d)=1}}\sum_{\substack{y\bmod c\\ady\equiv-m\bmod c}}\ue\Big(\frac{ny-ad}{c}\Big).\end{align*}
For each fixed $d\mid c$, from the condition $ady\equiv-m\bmod c$ it follows that $d\mid m$. In a similar way, we also conclude that $d\mid n$. We thus have
\begin{align*}
S(m,n;c)&=\sum_{d\mid (m,n,c)}\sum_{\substack{a\bmod{c/d}\\(a,c/d)=1}}\sum_{\substack{y\bmod c\\y\equiv-\overline{a}m/d\bmod{c/d}}}\ue\Big(\frac{ny-ad}{c}\Big)\\
&=\sum_{d\mid (m,n,c)}d\sum_{\substack{a\bmod{c/d}\\(a,c/d)=1}}\ue\Big(\frac{\overline{a}mn/d+ad}{c}\Big)\\
&=\sum_{d\mid (m,n,c)}dS(mnd^{-2},1;cd^{-1}).\end{align*}
This establishes \eqref{eq:KSidentity}.

\smallskip

\section{Comments}

The above argument also applies to the following general sum
\begin{align*}
\Xi_k(m,n;c)=\mathop{\sum\sum}_{\substack{x,y\bmod c\\xy\equiv k\bmod c}}\ue\Big(\frac{mx+ny}{c}\Big)
\end{align*}
with $k\in\bZ.$
In particular, we have $\Xi_1(m,n;c)=S(m,n;c).$ The identical arguments yield
\begin{align*}
\Xi_k(m,n;c)
&=\sum_{d|(m,n,c)}dS(mnd^{-2},k;cd^{-1}).
\end{align*}
On the other hand, the similar arguments (only using Lemma \ref{lm:orthogonality}) can also imply
\begin{align*}
\Xi_k(m,n;c)
&=\sum_{d|(m,k,c)}dS(mkd^{-2},n;cd^{-1}).
\end{align*}
Hence we see that $\Xi_k(m,n;c)$
is stable under permutations among $m,n,k.$ This is at least observed by Andersson \cite{An06}. This exponential sum $\Xi_k(m,n;c)$ appears in the work of Heath-Brown \cite{HB79} on the fourth moment of Riemann zeta functions and in Bykovsky, Kuznetsov and Vinogradov \cite{BKV90} when generalizing Kuznetsov's trace formula.

It is also interesting to see if the above arguments can be modified to derive a similar identity for the twisted sum
\begin{align*}
S_\chi(m,n;c)
&=\sideset{}{^*}\sum_{x\bmod c}\chi(x)\ue\Big(\frac{mx+n\overline{x}}{c}\Big).\end{align*}
Following similar arguments, we are also able to derive a Selberg type identity for $S_\chi(m,n;c)$ in the case that $\chi\bmod N$ and $(c/N,N)=(n,N)=1$ (as in \cite{An06}). However, new ideas are highly desired in general situations.

\smallskip

\bibliographystyle{plainnat}

\end{document}